\documentclass[fleq]{amsart}
\usepackage{amssymb,amsmath,latexsym,amsbsy,color,enumerate}
\numberwithin{equation}{section}

\newtheorem{theorem}{Theorem}

\newtheorem{lemma}{Lemma}

\newtheorem{proposition}{Proposition}
\newtheorem{corollary}{Corollary}

\newtheorem{conjecture}[theorem]{Conjecture}
\begin{document}
 \title{Local growth of pluri-subharmonic functions}
\author{Tuyen Trung Truong}
    \address{Department of mathematics, Indiana University, Bloomington IN 47405}
 \email{truongt@indiana.edu}
\thanks{}
    \date{\today}
    \begin{abstract}
We obtain two-bound estimates for the local growth of pluri-subharmonic functions in terms of Siciak and relative extremal functions. As applications, we
give simple new proofs of "Bernstein doubling inequality" and the main result in [Alexander Brudnyi,  Local inequalities for pluri-subharmonic functions,
Annals Math. 149 (1999), No. 2, pp. 511--533]. We propose a conjecture similar to the comparison theorem in [H. Alexander and B. A. Taylor,  Comparison
of two capacities in $\mathbb{C}^n$, Math. Z. 186 (1984), 407--417], whose validity allows to obtain bounds for the local growth of pluri-subharmonic
functions solely in term of the Siciak extremal functions.
 \end{abstract}
\maketitle
\section{Introduction}       
\label{SecIntroduction} Let $\Omega$ be an open subset of $\mathbb{C}^n$. The set of pluri-subharmonic functions on $\Omega$ is denoted as usual by
$PSH(\Omega )$. We are interested in obtaining bounds for the local growth of functions in $PSH(\Omega )$. Given two non-pluripolar sets
$A,E\subset\subset \Omega$, we define the function:
\begin{equation}
h_E(z):=\sup \{f(z)-\sup _{E}f:~f\in PSH(\Omega ),~\sup _{\Omega}f\leq 0,~\sup _{A}f\geq -1 \},\label{EquationFunctionH}
\end{equation}
where $z\in\Omega$. The problem is to obtain good estimates of the function $h_E(z)$ in terms of some intrinsic quantities of the set $E$, such as
(Lebesgue or Hausdorff) measures, or (logarithmic or relative) capacities. In this paper we will give some bounds of the function $h_E(z)$ by the later
quantities, via the Siciak and relative extremal functions. Let us recall the definitions of these extremal functions. The Siciak extremal function $V_E$
is defined as follows: For $z\in \mathbb{C}^n$
$$V_E(z)=\sup \{f(z):f\in \mathcal{L}(\mathbb{C}^n),~f|_E\leq 0\},$$
where $\mathcal{L}(\mathbb{C}^n)$ is the Lelong class
$$\mathcal{L}(\mathbb{C}^n)=\{f\in PSH(\mathbb{C}^n):~f(z)\leq \log^+ |z|+O(1)\}.$$
The relative extremal function $u_{E,\Omega }$ is defined as
$$u_{E,\Omega}(z)=\sup \{f(z):~f\in PSH(\Omega ),~f\leq 0,~ \sup _{E}f\leq -1\},$$
where $z\in\Omega$.

Our first result is
\begin{lemma}
i) We have
\begin{equation}
\frac{V_E(z)}{\sup _{\Omega}V_A}\leq h_E(z)\leq \frac{u_{E,\Omega}(z)+1}{|\sup _{A}u_{E,\Omega}|}. \label{EquationLemma1.1}
\end{equation}
ii) If $E$ is such that $u_{E,\Omega}$ is a continuous function then
\begin{eqnarray*}
h_E(z)= \frac{u_{E,\Omega}(z)+1}{|\sup _{A}u_{E,\Omega}|}.
\end{eqnarray*}
 \label{Lemma1}\end{lemma}

As some applications of Lemma \ref{Lemma1}, we will give simple new proofs to the main result in \cite{bru} and to the "Berstein doubling inequality".
The notation $B(x,\rho )$ (respectively $B_c(x,\rho )$) denotes the Euclidean ball with center $x$ and radius $\rho$ in $\mathbb{R}^n$ (respectively
$\mathbb{C}^n$). Let $r>1$ be a constant. Define $\mathcal{F}_r$ to be the set of functions $f\in PSH(B_c(0,r))$ satisfying
\begin{eqnarray*}
sup _{B_c(0,r)}f\leq 0,~sup _{B_c(0,1)}f\geq -1.
\end{eqnarray*}

\begin{theorem}
(Theorem 1.2 in \cite{bru}) Let the ball $B(x,t)$ satisfy $B(x,t)\subset B_c(x,at)\subset B_c(0,1)$, where $a>1$ is a fixed constant. There are constants
$c=c(a,r)$, $d=d(n)$ such that the inequality
\begin{equation}
\sup _{B(x,t)}f\leq c\log \frac{d|B(x,t)|}{|E|}+\sup _Ef, \label{EquationBrudnyiResult}\end{equation} holds for every $f\in \mathcal{F}_r$, and every
measurable set $E\subset B(x,t)$. (Here $|B(x,t)|$ and $|E|$ mean the Lebesgue measures of $B(x,t)$ and $E$, respectively, as subsets of $\mathbb{R}^n$.)
\label{TheoremBrudnyi}\end{theorem}
\begin{proposition}
(Proposition 2.5 in \cite{bru}) Let $f\in \mathcal{F}_r$ and $s\in [1,a]$, $a>1$. Suppose that $B_c(x,t)\subset B_c(x,at)\subset B_c(0,1)$. Then there is
a constant $c=c(r)$ such that
\begin{eqnarray*}
\sup _{B_c(x,st)}f\leq c\log s+\sup _{B_c(x,t)}f.
\end{eqnarray*}
 \label{PropositionBernstein}\end{proposition}

Let us remark that already in \cite{bru}, it was proved that when $n=1$, in the RHS of (\ref{EquationBrudnyiResult}) we can replace $|E|$ by the Siciak
capacity $C(E)$ of $E$. This suggests that for general $n$, we may obtain a similar result. We propose the following conjecture, whose validity allows
such an extension of Theorem \ref{TheoremBrudnyi} to the general cases when $E$ needs not to have positive Lebesgue measure.

\begin{conjecture} Let $A=B_c(0,1)$ and $\Omega =B_c(0,a)$. There exists a constant $C_{a,n}>0$ such that for all compact non-pluripolar set $E\subset A$ we have
\begin{equation}
|\sup _{A}u_{E,\Omega }|\sup _{\Omega}V_E\geq C_{a,n}.\label{EquationConjectureMain}
\end{equation}
 \label{ConjectureMain}\end{conjecture}

Let $\gamma =C(E)$ be the Siciak capacity of $E$, i.e.
$$\limsup _{s\rightarrow \infty}(\sup _{B_{c}(0,s)}V_E-\log s)=-\log \gamma .$$
The following is a corollary of conjecture \ref{ConjectureMain}.
\begin{corollary}
If conjecture \ref{ConjectureMain} is true, and if $\Omega =B_c(0,a)$, $A=B_c(0,1)$ then there exists $C_{a,n}>0$ such that for all compact
non-pluripolar set $E\subset B_c(0,1)$ we have:

\begin{equation} \frac{1}{C_{a,n}}\log \frac{1}{\gamma}\leq \sup _Ah_E\leq C_{a,n}\log \frac{n}{\gamma}.\label{EquationCorollary2.1}
\end{equation}
\label{Corollary2}\end{corollary} By Proposition \ref{PropositionBernstein}, as argued in \cite{bru} (see also the proof of Theorem \ref{TheoremBrudnyi}
in this paper), we can reduce proving (\ref{EquationBrudnyiResult}) to estimating
\begin{equation}
\sup _{B(0,1)}f-\sup _{E}f,\label{EquationCorollary2.2}
\end{equation}
where $f\in PSH(B_c(0,a))$, $\sup _{B_c(0,a)}f\leq 0$, $\sup _{B_c(0,1)}f\geq -1$ . Since the middle term of (\ref{EquationCorollary2.1}) is an upper
bound for the quantity in (\ref{EquationCorollary2.2}), Corollary \ref{Corollary2} may be viewed as an extension of Theorem \ref{TheoremBrudnyi}. Here
the set $E$ needs not to be a subset of $\mathbb{R}^n$  or to have positive ($\mathbb{R}^n$ or $\mathbb{C}^n$) Lebesgue measure.

Remark that conjecture \ref{ConjectureMain} is similar to the comparison theorem of Alexander-Taylor\cite{at}: There exists constants $c_n>0,~c_a>0$
(here $c_n$ depends only on $n$ and $c_a$ depends only on $a$) such that for all non-pluripolar set $E\subset A$ we have
\begin{equation}
\frac{c_n}{cap(E;\Omega )^{1/n}}\leq \sup _{A}V^*_E\leq \frac{c_a}{cap (E;\Omega )},\label{EquationAlexanderTaylorInequalities}
\end{equation}
where $cap(E;\Omega )$ is the relative capacity (for the definition, see for example \cite{at}). Note that the exponents of $cap(E;\Omega )$ in
(\ref{EquationAlexanderTaylorInequalities}) can not be improved. As explained in \cite{at}, the exponent $1/n$ in the LHS of
(\ref{EquationAlexanderTaylorInequalities}) occurs when $E$ is a ball, while the exponent $1$ in the RHS of (\ref{EquationAlexanderTaylorInequalities})
occurs when $E$ is a small polydisk. More generally, if $E=E_1\times \ldots \times E_n$ where $E_j\subset \mathbb{C}$, then in general the exponent may
be any number between $1/n$ and $1$. As will be shown later, in all these cases, conjecture \ref{ConjectureMain} holds. It is interesting to observe that
if $E$ is a ball of center $0$, then the LHS of (\ref{EquationConjectureMain}) is the constant $\log a$.

The rest of this paper is organized as follows. In Section 2, we prove Lemma \ref{Lemma1}, we prove Theorem \ref{TheoremBrudnyi} and Proposition
\ref{PropositionBernstein}. In Section 3, we verify conjecture \ref{ConjectureMain} in some cases, and prove Corollary \ref{Corollary2}.

\textbf{Acknowledgements.} The author would like to thank Professor Norman Levenberg for his generous help. The author also would like to thank Professor
Alexander Brudnyi for helpful comments.

\section{Proofs of Lemma \ref{Lemma1}, Proposition \ref{PropositionBernstein} and Theorem \ref{TheoremBrudnyi}}

Proof of Lemma \ref{Lemma1}

\begin{proof} i) Let $f\in PSH(\Omega )$ be such that $f\leq 0$, $\sup _Af\geq -1$. Define
$$\alpha :=\sup _Ef.$$
Then by the definition of $u_{E,\Omega }$ we have
$$f(x)\leq |\alpha |u_{E,\Omega }(x)=|\alpha |(u_{E,\Omega }(x)+1)+\alpha .$$
Hence
$$f(x)-\sup _{E}f=f(x)-\alpha \leq |\alpha |(u_{E,\Omega }(x)+1).$$
Now we estimate $|\alpha |$. We have
$$0\geq |\alpha |\sup _{A}u_{E,\Omega }\geq \sup _{A}f\geq -1.$$
Hence
$$|\alpha |\leq \frac{1}{|\sup _{A}u_{E,\Omega }|}.$$
Combining these inequalities we obtain
$$f(x)-\sup _{E}f\leq \frac{u_{E,\Omega}(x)+1}{|\sup _{A}u_{E,\Omega }|}.$$
Take supremum on over all such $f$, we obtain the RHS inequality of (\ref{EquationLemma1.1}).

Now we prove the LHS of (\ref{EquationLemma1.1}). Let $f\in \mathcal{L}(\mathbb{C}^n)$ be not a constant function with $\sup _{E}f= 0$. Consider the
function
$$g(z)=\frac{f(z)-\sup _{\Omega}f}{\sup _{\Omega}f-\sup _{A}f}.$$
Then $g\in PSH(\Omega )$, $\sup _{\Omega }g\leq 0$ and $\sup _Ag=-1$. Hence by definition of Siciak extremal function, we have
\begin{eqnarray*}
\frac{f(z)}{\sup _{\Omega}V_A}&\leq&
\frac{f(z)}{\sup _{\Omega}f-\sup _{A}f}\\
&=&g(z)-\sup _Eg\leq h_E(z).
\end{eqnarray*}
If we take supremum of the above inequality on over all such $f$ we obtain the LHS inequality of (\ref{EquationLemma1.1}).

ii) If $E$ is such that $u_{E,\Omega }$ is a continuous function then $u_{E,\Omega }$ itself is pluri-subharmonic in $\Omega$. Consider the function
$$g(z)=\frac{u_{E,\Omega }(z)}{|\sup _A u_{E,\Omega }|},$$
where $z\in \Omega$. Then $g\in PSH(\Omega )$, $\sup _{\Omega }g\leq 0$ and $\sup _Ag=-1$. Thus by definition of the $h_E$ we have
$$\frac{u_{E,\Omega}(z)+1}{|\sup _{A}u_{E,\Omega}|}=g(z)-\sup _{E}g\leq h_E(z).$$
\end{proof}

Proof of Proposition \ref{PropositionBernstein}:
\begin{proof}
In this case $\Omega =B_c(0,r)$, $A=B_c(0,1)$ and $E=B_c(x,t)$.

By Lemma \ref{Lemma1} we have
\begin{equation}
\sup _{B_c(x,st)}f\leq \frac{\sup _{B_c(x,st)}u_{B_c(x,t),B_c(0,r)}+1}{|\sup _{B_c(0,1)}u_{B_c(x,t),B_c(0,r)}|}+\sup
_{B_c(x,t)}f.\label{PropositionBernstein.1}
\end{equation}

By Proposition 5.3.3 in \cite{klimek} we have
\begin{eqnarray*}
\sup _{B_c(x,st)}u_{B_c(x,t),B_c(0,r)}+1\leq \frac{\sup _{B_c(x,st)}V_{B_c(x,t)}}{\inf _{\partial B_c(0,r)}V_{B_c(x,t)}}.
\end{eqnarray*}
Since $V_{B_c(x,t)}(z)=\log ^+(|z-x|/t)$, we obtain
\begin{eqnarray*}
\sup _{B_c(x,st)}u_{B_c(x,t),B_c(0,r)}+1\leq \frac{\log s}{\log ((r-1+t)/t)}.
\end{eqnarray*}

Now we estimate $|\sup _{B_c(0,1)}u_{B_c(x,t),B_c(0,r)}|$. Fix $z_0\in \partial B_c(0,1)$. We choose $l_{z_0}$ to be the complex line containing both
points $x$ and $z_0$. Then
\begin{equation}
|u_{B_c(x,t),B_c(0,r)}(z_0)|\geq |u_{B_c(x,t)\cap l_{z_0},B_c(0,r)\cap l_{z_0}}(z_0)|\geq |\sup _{B_c(0,1)\cap l_{z_0}}u_{B_c(x,t)\cap
l_{z_0},B_c(0,r)\cap l_{z_0}}|. \label{PropositionBernstein.2}\end{equation}

Now by the $1$-dimensional case of conjecture \ref{ConjectureMain}, which is known to be true (see for example \cite{at} or \cite{bru}, see also Section
4 in this paper), since $B_c(x,t)\cap l_{z_0}$ is a $1$-dimensional ball of radius $t$, there is a constant $C=C(r)$ depending only on $r$ such that
$$|\sup _{B_c(0,1)\cap l_{z_0}}u_{B_c(x,t)\cap l_{z_0},B_c(0,r)\cap
l_{z_0}}|\geq C/\sup _{B_c(0,r)\cap l_{z_0}}V_{B_c(x_0,t)\cap l_{z_0}}\geq C/\log ((r+1-t)/t).$$

Since $t\in [0,1]$, substituting all these inequalities into (\ref{PropositionBernstein.1}) we obtain
$$\sup _{B_c(x,st)}f\leq \frac{1}{C}\frac{\log ((r+1-t)/t)}{\log (r-1+t)/t}\log s+\sup
_{B_c(x,t)}f\leq C_1\log s+\sup _{B_c(x,t)}f,$$ where $C_1>0$ is a constant depending only on $r$.
\end{proof}

Proof of Theorem \ref{TheoremBrudnyi}
\begin{proof}
Using the "Bernstein doubling inequality" (Proposition \ref{PropositionBernstein}), as observed in \cite{bru} (page 523) it suffices to prove the
following equivalent statement. Let $a>1$ be a constant. Let $\mathcal{R}_a$ be the family of pluri-subharmonic functions on $B_c(0,a)$ satisfying the
conditions
$$\sup _{B_c(0,a)}f\leq 0,~\sup _{B_c(0,1)}f\geq -1.$$ Then for every measurable subset $E\subset B(0,1)$ of positive measure and every $f\in \mathcal{R}_a$
\begin{equation}
\sup _{B(0,1)}f\leq c\log \frac{d|B(0,1)|}{|E|}+\sup _{E}f.\label{EquationTheorem1.2}
\end{equation}
Here $d=d(n)$ and $c=c(a)$.

In this case $\Omega =B_c(0,a)$ and $A=B_c(0,1)$. Apply Lemma \ref{Lemma1} we get
\begin{equation}
\sup _{B(0,1)}f\leq \frac{\sup _{B(0,1)}u_{E,B_c(0,a)}+1}{|\sup _{B_c(0,1)}u_{E,B_c(0,a)}|}+\sup _{E}f.\label{EquationTheorem1.3}
\end{equation}

We divide the estimation of the first term in the RHS of (\ref{EquationTheorem1.3}) into several steps:

Step 1: $$|\sup _{B_c(0,1)}u_{E,B_c(0,a)}|\geq |\sup _{B_c(0,1)}u_{B(0,1),B_c(0,a)}|.|\sup _{B(0,1)}u_{E,B_c(0,a)}|.$$

Proof: Let $f$ be any function in $PSH(B_c(0,a) )$ with $\sup _{B_c(0,a) }f\leq 0$ and $\sup _{E}f\leq -1$. Define the function
$$g(z)=\frac{f(z)}{|\sup _{B(0,1)}u_{E,B_c(0,a)}|}.$$
Then $g\in PSH(B_c(0,a) )$, $\sup _{B_c(0,a) }g\leq 0$, and since $f(z)\leq u_{E,B_c(0,a)}(z)$ we have also $\sup _{B(0,1)}g\leq -1$. Thus by definition
of the relative extremal function
$$\frac{f(z)}{|\sup _{B(0,1)}u_{E,B_c(0,a)}|}=g(z)\leq u_{B(0,1),B_c(0,a)}(z),$$
for all $z\in \Omega$. Take supremum of the above inequality on over all such functions $f$, we obtain
$$\frac{u_{E,B_c(0,a)}(z)}{|\sup _{B(0,1)}u_{E,B_c(0,a)}|}\leq u_{B(0,1),B_c(0,a)}(z).$$
Form this we obtain the claim of Step 1.

Step 2: Apply Step 1 to (\ref{EquationTheorem1.3}), for any $f\in \mathcal{R}_a$ we have
\begin{equation}
\sup _{B(0,1)}f\leq C_1\frac{\sup _{B(0,1)}u_{E,B_c(0,a)}+1}{|\sup _{B(0,1)}u_{E,B_c(0,a)}|}+\sup _{E}f,\label{EquationTheorem1.4}
\end{equation}
where
$$C_1=\frac{1}{|\sup _{B_c(0,1)}u_{B(0,1),B_c(0,a)}|},$$
depends only on $a$.

Step 3: Let $x_0$ be any point in $B(0,1)$. Then by Lemma 3 of \cite{bru-ganzburg}, there exists a ray $l_0$ such that

\begin{equation}
\frac{mes _1(B(0,1)\cap l_0)}{mes _1(E\cap l_0)}\leq \frac{n|B(0,1)|}{|E|}.\label{EquationTheorem1.5}
\end{equation}

Let $l_0'$ be the one-dimensional affine complex line containing $l_0$. Using the properties of extremal functions in one-dimensional and
(\ref{EquationTheorem1.5}), we obtain
\begin{eqnarray*}
\frac{\sup _{B(0,1)}u_{E,B_c(0,a)}+1}{|\sup _{B(0,1)}u_{E,B_c(0,a)}|}&=&\sup _{z_0\in
B(0,1)}\frac{u_{E,B_c(0,a)}(z_0)+1}{|u_{E,B_c(0,a)}(z_0)|}\\
&\leq& \sup _{z_0\in B(0,1)}\frac{u_{E\cap l_0',B_c(0,a)\cap l_0'}(z_0)+1}{|u_{E\cap
l_0',B_c(0,a)\cap l_0'}(z_0)|}\\
&\leq&\sup _{z_0\in B(0,1)}\frac{V_{E\cap l_0'}(z_0)}{|u_{E\cap
l_0',B_c(0,a)\cap l_0'}(z_0)|\inf _{\partial (B_c(0,a)\cap l_0')}V_{E\cap l_0'}} \\
 &\leq &C_2\log \frac{4~ mes_1(B(0,1)\cap l_0)|}{mes_1(E\cap l_0)} \leq C_2\log \frac{4n|B(0,1)|}{|E|},
\end{eqnarray*}
for some constant $C_2>0$ depending only on $a$. This inequality together with (\ref{EquationTheorem1.4}) complete the proof of Theorem
\ref{TheoremBrudnyi}.
\end{proof}

\section{Verification of conjecture \ref{ConjectureMain} in some cases}
Throughout this section $\Omega =B_c(0,a)$, $A=B_c(0,1)$ and $E$ is a compact subset of $A$.

We need the following results

Claim 1:  \begin{equation} \log \frac{1}{\gamma}\leq \sup _AV_E\leq 2e^2n\log \frac{n}{\gamma}.\label{EquationTaylorInequality}
\end{equation}
\begin{proof}
The LHS of (\ref{EquationTaylorInequality}) follows easily from the following two facts:

i) If $s\geq t>0$ then
$$\sup _{B_{c}(0,s)}V_E-\log s\leq \sup _{B_{c}(0,t)}V_E-\log t.$$

ii) $$\limsup _{s\rightarrow \infty}\sup _{B_{c}(0,s)}V_E-\log s=-\log \gamma .$$

The proof of the RHS of (\ref{EquationTaylorInequality}) is similar to the proof of formula (1.2) in \cite{tuyen}: we use Taylor's inequality (see
\cite{taylor}) applied to estimate the integration of $V_E^*$ on the sphere $|z|=n$, and the Harnack inequality for positive PSH functions.
\end{proof}

Claim 2:  $$\sup _Au_{E,\Omega }+1\leq 2\frac{\sup _AV_E}{\sup _{\Omega }V_E}.$$.
\begin{proof}
Define $M=\sup _{\Omega }V_E$. For a function $u$, let $u^*$ be the upper-semicontinuous regularization of $u$. Then it is well-known that the function
$V_E^*$ is in the Lelong class $\mathcal{L}(\mathbb{C}^n)$. Consider the following function

$$V(z)=(\sup _{B_c(0,|z|)}V_E^*)^*.$$
Then $V(z)$ is also in the Lelong class $\mathcal{L}(\mathbb{C}^n)$.

Fix a function $f\in PSH(\Omega )$ with $\sup _{\Omega }f\leq 0$, $\sup _{E}f\leq -1$. Define
$$u(z)=\left \{\begin{array}{ll}\max\{M(f(z)+1),V(z)\},&z\in \Omega \\V(z),&z\in \mathbb{C}^n\backslash \Omega .\end{array}\right .$$
Then $u(z)$ is in the Lelong class $\mathcal{L}(\mathbb{C}^n)$. Hence
$$u(z)\leq V_E(z)+\sup _Eu.$$
Now we estimate $\sup _Eu$. Since $E\subset A=B_c(0,1)$, we have:
$$\sup _{E}u=\sup _EV\leq \sup _AV=\sup _AV_E^*=\sup _AV_E.$$
In particular
$$M(f(z)+1)\leq V_E(z)+\sup _AV_E.$$
Take supremum on over all such $f$, we obtain
$$M(u_{E,\Omega }+1)\leq V_E(z)+\sup _AV_E.$$
Thus
$$\sup _Au_{E,\Omega }+1\leq 2\frac{\sup _AV_E}{\sup _{\Omega }V_E}.$$
\end{proof}

We verify conjecture \ref{ConjectureMain} in the following four cases:

Case 1: $n=1$. In this case Conjecture \ref{ConjectureMain} is just the Alexander-Taylor inequality (\ref{EquationAlexanderTaylorInequalities}), using
the equivalence between $cap(E;\Omega )$ and $|\sup _{A}u_{E,\Omega}|$ (see \cite{at}).

Case 2: $E=\prod _{j=1}^nD_j$ is a polydisk, where $D_j$ is a disk in $\mathbb{C}$. In this case the Siciak capacity $\gamma =Cap(E)$ of $E$ is the
smallest radius of the disks $D_j$'s. The same argument as that of the proof of Proposition \ref{PropositionBernstein}, together with
(\ref{EquationTaylorInequality}),  proves conjecture \ref{ConjectureMain} in this case.

Case $3$: $E\subset B_c(z_0,\gamma ^{\tau _n})$ where $\gamma =Cap(E)$ is the Siciak capacity of $E$, and
$$\tau _n=1-\frac{1}{8e^2n}.$$
Without loss of generality (using the automorphism of $\Omega$ translating $z_0$ to the origin $0\in \mathbb{C}^n$), we may assume that $z_0=0$. It
suffices to prove Conjecture \ref{ConjectureMain} when $\gamma$ is small enough.

The proof of Claim 2 and (\ref{EquationTaylorInequality}) gives
\begin{eqnarray*}
\sup _{B_c(0,\gamma ^{\tau _n})}u_{E,\Omega }&\leq& 2\frac{\sup _{B_c(0,\gamma ^{\tau _n})}V_{E}}{\sup _{\Omega }V_E}-1\leq 4e^2n(1-\tau _n)\frac{-\log
\gamma }{\log a-\log \gamma}-1.
\end{eqnarray*}
Hence when $\gamma$ is small enough we have
$$\sup _{B_c(0,\gamma ^{\tau _n})}u_{E,\Omega }\leq -\frac{1}{3}.$$
Then it follows that
$$|\sup _Au_{E,\Omega }|\geq \frac{1}{3}|\sup _Au_{B_c(0,\gamma ^{\tau _n}),\Omega}|.$$
This inequality, together with the LHS of (\ref{EquationTaylorInequality}) completes the proof of Conjecture \ref{ConjectureMain} for Case 3.

Remark: A similar constraint was used in \cite{tuyen} (see Lemma 1 in \cite{tuyen}) when exploring sets non-thin at $\infty$ in $\mathbb{C}^n$.

Case 4: $E=E_1\times \ldots \times E_n$, where $E_j\subset \mathbb{C}$ are compact non-pluripolar, and $a>\sqrt{n}$. In this case, there exists $r>1$
such that $A=B_c(0,1)\subset B=D(0,r)\times D(0,r)\ldots \times D(0,r)\subset \Omega =B_c(0,a)$, where $D(0,r)\subset \mathbb{C} $ is the one-dimensional
disk. Then
$$|\sup _{A}u_{E,\Omega }|\geq |\sup _A u_{E,B}|=|\sup _A u^*_{E,B}|.$$
We also have
$$\sup _{A}V_{E}= \sup _A V^*_{E}.$$
Using the product property of the function $u^*_{E,B}$ and $V_E^*$ (see for example \cite{edi-pol} and \cite{bloki}), Case 4 is reduced to Case 1 above.

Proof of Corollary \ref{Corollary2}: From Lemma \ref{Lemma1} and the arguments above, Corollary \ref{Corollary2} follows easily.

\end{document}